\documentclass[12pt]{article}

\usepackage{amsmath}  
\allowdisplaybreaks 
\usepackage{amssymb, amsthm, mathtools}

\usepackage{authblk}

\usepackage{float}

\usepackage{booktabs}
\usepackage[colorlinks=true, urlcolor=blue, linkcolor=blue]{hyperref}

\usepackage[ruled,vlined,linesnumbered]{algorithm2e}
\SetKwInput{KwIn}{Input}
\SetKwInput{KwOut}{Output}
\SetKwComment{tcp}{\% }{}

\newcommand{\proj}{\operatorname{proj}}

\usepackage{geometry}
\title{\textbf{A Generalized $(k,m)$ Heron Problem: Optimality Conditions and Algorithm}}

\author[1]{Triloki Nath\thanks{Department of Mathematics and Statistics, Deen Dayal Upadhyaya Gorakhpur University, Gorakhpur, India. 
		Email: \texttt{tnverma07@gmail.com}}}
\author[2]{Manohar Choudhary\thanks{Department of Mathematics and Statistics, Dr. Hari Singh Gour Vishwavidyalaya, Sagar, India. 
		Email: \texttt{manoharfbg@gmail.com}}}
\author[2]{Ram K. Pandey\thanks{Department of Mathematics and Statistics, Dr. Hari Singh Gour Vishwavidyalaya, Sagar, India. 
		Email: \texttt{pandeywavelet@gmail.com}}}
\affil[1]{Department of Mathematics and Statistics, Deen Dayal Upadhyaya Gorakhpur University, Gorakhpur, India}
\affil[2]{Department of Mathematics and Statistics, Dr. Hari Singh Gour Vishwavidyalaya, Sagar, India}

\date{}

\theoremstyle{plain}
\newtheorem{theorem}{Theorem}[section]
\newtheorem{lemma}[theorem]{Lemma}
\newtheorem{proposition}[theorem]{Proposition}

\theoremstyle{definition}
\newtheorem{definition}[theorem]{Definition}
\newtheorem{remark}[theorem]{Remark}
\newtheorem{example}[theorem]{Example}

\usepackage{rotating}  
\usepackage{booktabs}  
\usepackage{hyperref}

\begin{document}
	
\maketitle

\begin{abstract}
	This paper presents a new extension of the classical Heron problem, termed as the \emph{generalized $(k,m)$-Heron problem}, which seeks an optimal configuration among $k$ feasible and $m$ target non-empty closed convex sets in $\mathbb{R}^n$. The problem is formulated as to find a point in each set that minimize pairwise distances from the point in each $k-$feasible sets to the point in $m-$target sets. This leads to a convex optimization framework that generalizes several well-known geometric distance problems. Using tools from convex analysis, we establish fundamental results on the existence, uniqueness, and first-order optimality conditions through subdifferential calculus and normal cone theory. Building on these insights, a \emph{Projected Subgradient Algorithm (PSA)} is proposed for numerical solution, and its convergence is rigorously proved under a diminishing step-size rule. Numerical experiments in $\mathbb{R}^2$ and $\mathbb{R}^3$ illustrate the algorithm’s stability, geometric accuracy, and computational efficiency. Overall, this work provides a comprehensive analytical and algorithmic framework for multi-set geometric optimization with promising implications for location science, robotics, and computational geometry.
\end{abstract}

\noindent\textbf{Keywords:} Generalized Heron problem, Convexity, Projected subgradient algorithm, Distance Minimization.\\
\noindent\textbf{MSC Classifications:}  90B85, 51M04,

\section{Introduction}

Geometric optimization problems play a fundamental role in mathematics, engineering, computer science, and location science, where the objective is to determine the location or arrangement of points and sets that minimize or maximize certain geometric or distance-based criteria. Among the earliest and most celebrated examples is the \emph{Heron problem}, attributed to \textit{Heron of Alexandria (10–75 AD)}, which can be stated as follows:

\begin{quote}
	\textit{Given two points $A$ and $B$ lying on the same side of a straight line $L$ in the plane, find a point $C \in L$ such that the total distance $|CA| + |CB|$ is minimized.}
\end{quote}

The classical Heron problem exemplifies a \emph{distance minimization} task constrained to a geometric set. Over the centuries, this problem has inspired numerous extensions and generalizations, leading to rich developments in \emph{convex analysis, variational principles, and computational geometry}.

Recall that a set $C \subseteq \mathbb{R}^n$ is \emph{convex} if it contains all line segments joining any two of its points.  
Formally, $C$ is convex if 
$\lambda x + (1 - \lambda)y \in C,$  $\forall\, x, y \in C$ and  $0 \leq \lambda \leq 1.$

In Convex setting, the first comprehensive generalization of the classical Heron problem was introduced by \textbf{Mordukhovich, Nam, and Salinas}~\cite{mordukhovich2011Variational,mordukhovich2012convex}, who reformulated the problem within the framework of modern convex analysis. In their \emph{generalized Heron problem}, the two fixed points are replaced by finitely many nonempty, closed, and convex sets $\{C_1, C_2, \dots, C_m\} \subset \mathbb{R}^n$, while the constraint line is generalized to a nonempty, closed, and convex feasible set $S \subset \mathbb{R}^n$. The objective is to determine a point $x \in S$ that minimizes the total Euclidean distance to all target sets:
\begin{equation}\label{eq:genHeronBasic}
	\min_{x \in S} \; D(x) = \sum_{j=1}^m d(x, C_j),
\end{equation}
where $d(x, C_j)$ denotes the Euclidean distance from $x$ to the set $C_j$. 

This formulation casts the classical Heron problem into a convex optimization framework, where tools from subdifferential calculus yield necessary and sufficient optimality conditions. Moreover, the authors proposed a convergent \emph{projected subgradient algorithm (PSA)} for its numerical solution.

Subsequent progress was made by \textbf{Chi and Lange}~\cite{chi2014mm}, who revisited the generalized Heron problem from an algorithmic perspective. They employed the \emph{majorization–minimization (MM) principle} to develop fast and efficient iterative schemes for the Euclidean version of the problem, significantly improving the computational performance compared to subgradient-based methods. Their work, along with earlier contributions on related geometric problems such as the \emph{Fermat–Weber} problem~ \cite{kuhn1973note,MordukhovichNamFermat2011,weiszfeld1937fermat,WeiszfeldPlastria2009}  and the \emph{Steiner} problem~\cite{courant1961mathematics,gueron2002fermat}, illustrates the close connection among these classical formulations and the modern generalized Heron framework. Collectively, these contributions reveal a rich interplay between geometric modeling and convex optimization, motivating further generalizations of Heron-type problems to more complex and multi-set frameworks.

\subsection{The generalized $(k,m)$  Heron problem}

Motivated by these developments, we propose a \emph{further generalization} that extends the classical and generalized Heron problems to interactions between multiple feasible and multiple target convex sets in $\mathbb{R}^n$.

Let $S_1, S_2, \dots, S_k \subset \mathbb{R}^n$ denote nonempty, closed, and convex feasible sets, and let $C_1, C_2, \dots, C_m \subset \mathbb{R}^n$ denote nonempty, closed, and convex target sets.  
From each feasible set $S_i$, we select a point $x_i \in S_i$, and from each target set $C_j$, we select a point $y_j \in C_j$.  
The total \emph{interaction cost} is defined as the sum of pairwise Euclidean distances between these feasible and target points:
\begin{equation}\label{eq:objective}
	F(x, y) = \sum_{i=1}^k \sum_{j=1}^m \|x_i - y_j\|.
\end{equation}
The \emph{generalized $(k,m)$ Heron problem} is then formulated as the convex optimization problem:
\begin{equation}\label{eq:(k,m)HeronProblem}
	\begin{aligned}
		\min~  & F(x, y) = \sum_{i=1}^k \sum_{j=1}^m \|x_i - y_j\|, \\
		\text{s.t.} \quad & x_i \in S_i, \quad i=1,\dots,k, \\
		& y_j \in C_j, \quad j=1,\dots,m.
	\end{aligned}
\end{equation}

\begin{figure}[H]
	\centering
	\includegraphics[width=0.8\textwidth]{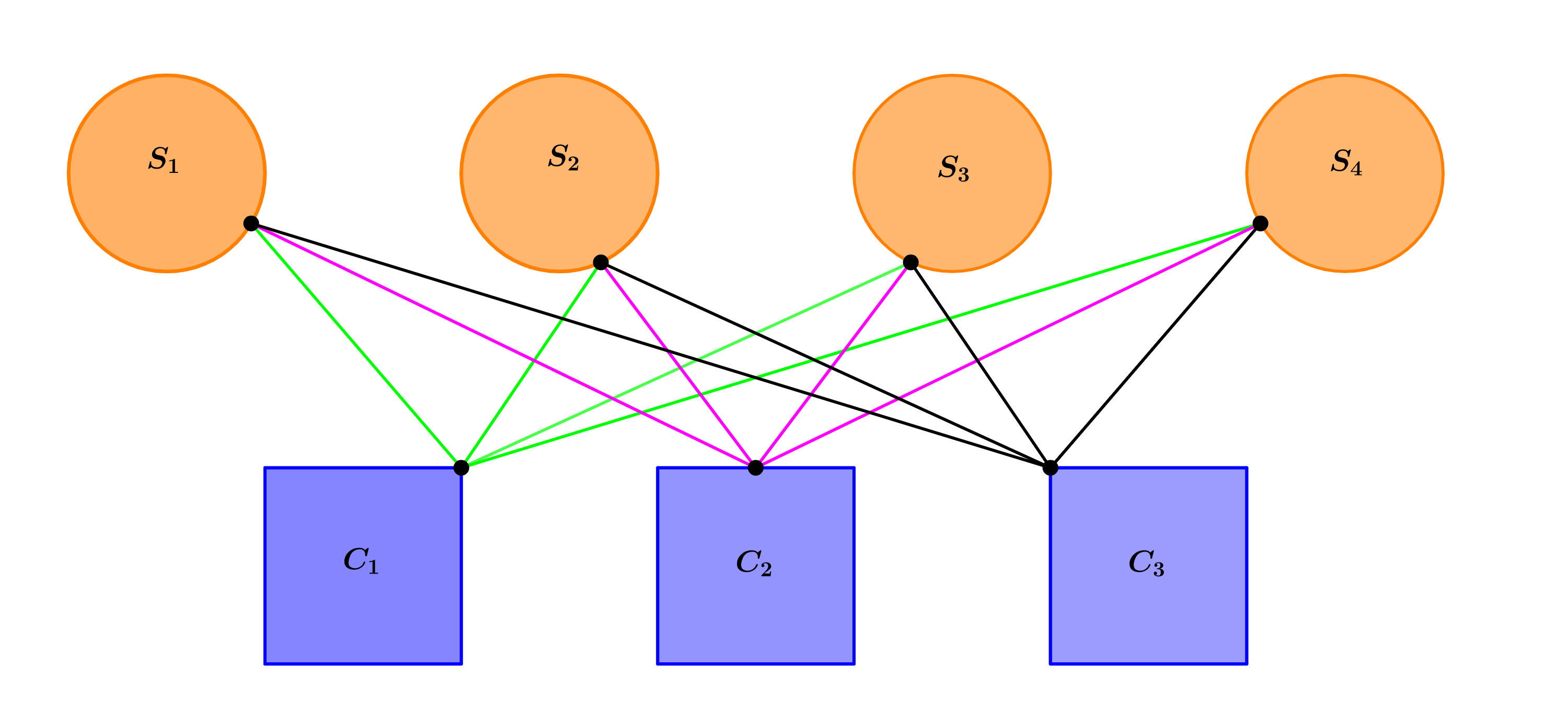}
	\caption{Illustration of generalized $(k,m)$ Heron problem.}
	\label{fig:kmHeronIllustration}
\end{figure}

\subsection{Connection to the generalized Heron problem}

The generalized $(k,m)$ Heron problem naturally reduces to the standard generalized Heron formulation when $k = 1$, that is, when the feasible region consists of a single closed and convex set $S_1 \subset \mathbb{R}^n$. In this case, problem~\eqref{eq:(k,m)HeronProblem} becomes
\[
\min_{x \in S_1,\, y_j \in C_j,\, j = 1, \dots, m} \; \sum_{j=1}^m \|x - y_j\|.
\]

For any fixed $x \in S_1$, the inner minimization over each $y_j \in C_j$ simplifies to
\[
\min_{y_j \in C_j} \|x - y_j\| = d(x, C_j),
\]
where the minimum is attained at the Euclidean projection $y_j^* = \proj_{C_j}(x)$.  
Since each $C_j$ is nonempty, closed, and convex, such a projection exists and is unique.  
Substituting these optimal points into the objective function reduces the problem to the well-known \emph{generalized Heron problem}~\cite{mordukhovich2012convex}:
\begin{equation}\label{eq:genHeron}
	\min_{x \in S_1} \; \sum_{j=1}^m d(x, C_j),
\end{equation}
which seeks a point $x \in S_1$ minimizing the sum of its Euclidean distances to the target sets $\{C_j\}_{j=1}^m$.

Hence, the proposed \emph{$(k,m)$-Heron problem} includes the classical and generalized Heron formulations as special cases, while introducing a broader framework that captures interactions among multiple feasible and target regions. This extension establishes a natural connection between multi-set geometric optimization and modern convex analysis.

\subsection{Contributions of this work}

Building upon the generalized $(k,m)$-Heron framework introduced above, this paper advances the study of the problem through both theoretical analysis and computational development.  
Section~\ref{sec:preliminaries} presents the essential \emph{preliminaries} from convex analysis to ensure the paper is self-contained and accessible.  
Section~\ref{optimality} establishes the \emph{existence}, \emph{uniqueness}, and \emph{characterization} of optimal solutions and derives the \emph{necessary and sufficient optimality conditions}.    
In Section~\ref{numerical-algorithm}, a \emph{Projected Subgradient Algorithm (PSA)} is developed, and its convergence is proved under a diminishing step-size rule $\alpha_t = 1/t$.  
Finally, Section~\ref{sec:numerical-example} presents two numerical examples in $\mathbb{R}^2$ and $\mathbb{R}^3$, illustrating the convergence behavior of the algorithm and supporting the theoretical results.

\section{Preliminaries}\label{sec:preliminaries}

This section recalls several fundamental concepts and results from convex analysis that form the mathematical foundation of the generalized $(k,m)$-Heron problem.  
Key properties of convex sets and convex functions---such as the fact that every local minimum of a convex function is also global, and that the distance function to a convex set is convex---are central to both the theoretical analysis and algorithmic development.   
For a comprehensive background on convex analysis, the reader is referred to the standard texts, e.g. ~\cite{Be14,boyd2004convex,HiLe96,MoNa23,Ro70}.

\begin{definition}
	For a function \( f : \mathbb{R}^n \to (-\infty, +\infty] \), the \emph{effective domain} is
	\[
	\mathrm{dom} f := \{ x \in \mathbb{R}^n \mid f(x) < +\infty \}.
	\]
\end{definition}

\begin{definition}
	A function \( f : \mathbb{R}^n \to (-\infty, +\infty] \) is \emph{convex} if \(\mathrm{dom} f\) is convex and
	\[
	f(\lambda x + (1 - \lambda) y) \le \lambda f(x) + (1 - \lambda) f(y)
	\]
	for all \(x, y \in \mathrm{dom} f\) and \( \lambda \in [0,1] \).  
	It is \emph{strictly convex} if the inequality is strict whenever \(x \ne y\) and \( \lambda \in (0,1) \).
\end{definition}
The following is the notion of a set associated with a function, which plays an important role in convex analysis.
\begin{definition}
	The \emph{epigraph} of \( f : \mathbb{R}^n \to \mathbb{R} \) is
	\[
	\mathrm{epi}\, f := \{ (x, t) \in \mathbb{R}^n \times \mathbb{R} \mid f(x) \le t \}.
	\]
	It is easy to verify that a function \(f\) is convex if and only if \(\mathrm{epi}\, f\) is a convex set.
\end{definition}

\begin{definition}
	For a nonempty set $C \subset \mathbb{R}^n$, the \emph{distance function} is
	\[
	d_C(x) := \inf_{y \in C} \|x-y\|.
	\]
	It is easy to see that if $C$ is convex, then $d_C$ is convex and Lipschitz continuous.
\end{definition}

\begin{definition}
	If $C \subset \mathbb{R}^n$ is nonempty, then \emph{projection} of $x \in \mathbb{R}^n$ onto $C$ is
	\[
	\proj_C(x) := \arg\min_{y \in C} \|x-y\|.
	\]
	For a nonempty closed convex set $C$, the projection is unique, and Moreover, $d_C(x) = \|x-\proj_C(x)\|$.
\end{definition}

Sometimes, auxiliary extended real-valued functions play a pivotal role in many problems.  
The following is one such function, called the \emph{indicator function}, which converts a constrained problem into an unconstrained one.

\begin{definition} 
	For a set $C \subset \mathbb{R}^n$, the \emph{indicator function} is
	\[
	\delta_C(x) := 
	\begin{cases}
		0, & x \in C, \\
		+\infty, & x \notin C.
	\end{cases}
	\]
\end{definition}

\begin{definition}  
	Let \( C \subset \mathbb{R}^n \) be convex, and let \( x \in C \). The \emph{normal cone} to \( C \) at \( x \) is
	\[
	N_C(x) := \{ v \in \mathbb{R}^n \mid \langle v, y - x \rangle \le 0 \ \ \forall y \in C \},
	\]
	where \(\langle \cdot, \cdot \rangle\) denotes the standard Euclidean inner product.  
	Geometrically, \( N_C(x) \) consists of all vectors that form an obtuse angle with every feasible direction from \(x\) within \(C\).
\end{definition}

The following result provides a nice tool to compute normal cone to cartessian product of convex sets, and it is necessary what we follows.

\begin{proposition}\cite[Page 59, Proposition 2.39]{DhDu11} 	\label{Normalcone_Product}
	Consider two closed convex sets $C_i \subset \mathbb{R}^{n_i}$, for $i = 1, 2$. 
	Let $a_i \in C_i$, for $i = 1, 2$. Then
	\begin{equation}
		N_{C_1 \times C_2}\bigl((a_1, a_2)\bigr)=  	N_{C_1}(a_1) 
		\times 	N_{C_2}(a_2).
	\end{equation}
	
\end{proposition}

Since the distance function is not differentiable when \(x \in C\), we require tools from non-smooth  analysis. The central object is the subgradient, which generalizes the gradient.

\begin{definition}
	Let $f:\mathbb{R}^n \to (-\infty,+\infty]$ be convex. A vector $v \in \mathbb{R}^n$ is a \emph{subgradient} of $f$ at $x \in \mathrm{dom}\,f$ if
	\[
	f(y) \ge f(x) + \langle v, y-x \rangle \quad \forall y \in \mathbb{R}^n.
	\]
	The set of all subgradients is called the \emph{subdifferential} of $f$ at $x$, denoted $\partial f(x)$.
\end{definition}

For the indicator function \(\delta_C\) of a convex set \( C \), the subdifferential at any \(x \in C\) coincides with the normal cone to \(C\) at \(x\)~\cite[Example~3.7, p.~85]{MoNa23}:
\begin{equation}\label{eq:indicator-subgrad}
	\partial \delta_C(x) = N_C(x), \quad x \in C.
\end{equation}

Investigating the GHP, the subdifferential of the distance function plays a crucial role. For algorithmic purposes, the computation of subdifferentials is essential

\begin{proposition}[Subdifferential of the distance function]
	\label{subgradientformula}
	Let \( C \subset \mathbb{R}^n \) be nonempty, closed, and convex. Then for any \(x \in \mathbb{R}^n\),
	\[
	\partial d_C(x) =
	\begin{cases}
		\left\{ \dfrac{x - \proj_C(x)}{\|x - \proj_C(x)\|} \right\}, & \text{if } x \notin C, \\[10pt]
		N_C(x) \cap \mathbb{B}, & \text{if } x \in C,
	\end{cases}
	\]
	where \(\mathbb{B} := \{ z \in \mathbb{R}^n \mid \|z\| \le 1 \}\) is the closed unit ball.  
	A proof can be found in~\cite[Example~3.3, p.~259]{HiLe96}.
\end{proposition}

In the special case where \( C = \{y\} \) is a singleton, we have \(\proj_C(x) = \{y\}\) and $d_C(x) = \|x - y\|.$ Then, for \(x \ne y\),
\begin{equation}
	\label{subgradientformula_forsingleton}
	\partial d_C(x) = \left\{ \frac{x - y}{\|x - y\|} \right\}.
\end{equation}

\begin{theorem}[Affine Composition Rule, {\cite[Theorem 4.2.1, p.263]{HiLe96}}]\label{affinetheorem}
	Let  $T: \mathbb{R}^n \to \mathbb{R}^m$ 
	be an affine mapping defined by $T(x) \;=\; T_0 \,x \;+\; c,$ where  	$T_0: \mathbb{R}^n \to \mathbb{R}^m$
	is a linear operator and  \(c \in \mathbb{R}^m\).  
	Let \(\psi: \mathbb{R}^m \to \mathbb{R}\) be a finite convex function. Then, for every \(x \in \mathbb{R}^n\),
	\[
	\partial\bigl(\psi\circ T\bigr)(x)
	\;=\;
	T_0^*\,\partial \psi\bigl(T(x)\bigr),
	\]
	where $T_0^*: \mathbb{R}^m \to \mathbb{R}^n$
	denotes the adjoint operator of \(T_0\). In other words, \(T_0^*\) is the unique linear map satisfying
	\[
	\langle T_0\,x,\; v \rangle
	\;=\;
	\langle x,\; T_0^*\,v \rangle
	\quad
	\forall\, x \in \mathbb{R}^n,\; v \in \mathbb{R}^m.
	\]
\end{theorem}

Following the discussion on \cite[p.\ 263]{HiLe96}, Let $f: \underbrace{\mathbb{R}^n \times \cdots \times \mathbb{R}^n}_{m\text{ times}}	\;\longrightarrow\; \mathbb{R}$
be a convex function in \(m\) blocks of variables, i.e.\ we write 
\((x_1, x_2, \dots, x_m)\) with each \(x_i \in \mathbb{R}^n\).  
Fix any \((x_2,\dots,x_m) \in \mathbb{R}^{n(m-1)}\), and define the affine map $ T: \mathbb{R}^n \to \mathbb{R}^{mn},\quad T(x_1) = (x_1, x_2, \dots, x_m)$ 
The linear part of \( T \) is \( T_0(x_1) = (x_1, 0, \dots, 0) \), and the corresponding adjoint is \( T_0^*(v_1, v_2, \dots, v_m) = v_1 \). Hence \( f \circ T \) is precisely the function \( x_1 \mapsto f(x_1, x_2, \dots, x_m) \).

\medskip

For each \(1 \le i \le m\), define the partial function
\[
f_i: \mathbb{R}^n \to \mathbb{R},
\quad
f_i(x_i) 
\;=\;
f\bigl(x_1,\dots,x_{i-1},\;x_i,\;x_{i+1},\dots,x_m\bigr),
\]
where the other \(x_j\) (with \(j \neq i\)) are held fixed. Applying~\ref{affinetheorem} to each such affine composition gives:
\[
\partial f_i(x_i)
\;=\;
\Bigl\{
v_i \in \mathbb{R}^n 
\;\Bigm|\;
\exists\,v_j \,(j\neq i)\text{ with }(v_1,\dots,v_m)\in \partial f(x_1,\dots,x_m)
\Bigr\}.
\]
Thus, \(\partial f_i(x_i)\) is the \emph{projection} of \(\partial f(x_1,\dots,x_m)\) onto the \(i\)th block. Consequently,
\begin{equation}\label{subset_inclusion}
	\partial f(x_1,\dots,x_m)
	\;\subseteq\;
	\partial f_1(x_1)
	\;\times\;
	\partial f_2(x_2)
	\;\times\;
	\cdots
	\;\times\;
	\partial f_m(x_m).
\end{equation}

\begin{remark}\label{remark:subgradient_equality}
	In general, the inclusion above in \ref{subset_inclusion}  is strict. However, if for each \(i\) the partial function \(f_i\) is differentiable at \(x_i\), then \(\partial f_i(x_i)\) reduces to a singleton \(\{\nabla f_i(x_i)\}\). Consequently, if all \(f_i\) are differentiable at the respective \(x_i\) for \(i = 1, \dots, m\), the subdifferential becomes as:
	\begin{equation}\label{sub_cartessian}
		\partial f(x_1,\dots,x_m)
		\;=\;
		\Bigl\{
		\bigl(\nabla f_1(x_1),\,\nabla f_2(x_2),\,\dots,\,\nabla f_m(x_m)\bigr)
		\Bigr\}.
	\end{equation}
\end{remark}

\begin{example}[Subdifferential of the pairwise distance]
	\label{ex:subdiff_distance}
	Let \( D(a_1, \dots, a_m) = \|a_1 - a_2\| \) denote the pairwise distance between \( a_1 \) and \( a_2 \) in \( \mathbb{R}^n \), $a_1 \ne a_2$ viewed as a function of the tuple \( (a_1, \dots, a_m) \in \mathbb{R}^{nm}\). The subdifferential of \( D \) with respect to the full tuple is given by the ordered tuple of partial subdifferentials:
	\[
	\partial D(a_1, \dots, a_m) = \bigl( \partial D_1(a_1),\, \partial D_2(a_2),\, \partial D_3(a_3),\, \dots,\, \partial D_m(a_m) \bigr),
	\]
	where \( \partial D_k(a_k) \) denotes the projection of \( \partial D(a_1, \dots, a_m) \) onto the \( k \)-th coordinate.
	
	By~\ref{remark:subgradient_equality} and the preceding discussion, we obtain:
	\[
	\partial D_1(a_1) = \frac{a_1 - a_2}{\|a_1 - a_2\|}, \qquad
	\partial D_2(a_2) = \frac{a_2 - a_1}{\|a_1 - a_2\|}, \qquad
	\partial D_k(a_k) = \mathbf{0} \quad \text{for all } k \ge 3.
	\]
	Applying the Cartesian product rule for subdifferentials see \eqref{sub_cartessian}, we conclude that
	\[
	\partial D(a_1, \dots, a_m) =
	\biggl(
	\frac{a_1 - a_2}{\|a_1 - a_2\|},\;
	\frac{a_2 - a_1}{\|a_1 - a_2\|},\;
	\mathbf{0}, \dots, \mathbf{0}
	\biggr).
	\]
	Thus, only the first two components contribute to the subdifferential while the remaining \( m - 2 \) components are zero.
\end{example}

Following is the most fundamental result regarding subdifferentials calculus and is widely used as sum rule.

\begin{theorem}[Moreau--Rockafellar Theorem {\cite[Corollary~3.21, p.~93]{MoNa23}}]
	\label{sumofsubgradient}
	Let \(\psi_i : \mathbb{R}^n \to (-\infty, +\infty]\), \(i=1,\dots,k\), be closed convex functions.  
	Suppose there exists \(\bar{x} \in \bigcap_{i=1}^k \mathrm{dom} \psi_i\) at which all but at most one of the functions are continuous.  
	Then for all \(x \in \bigcap_{i=1}^k \mathrm{dom} \psi_i\),
	\[
	\partial \left( \sum_{i=1}^k \psi_i \right)(x)
	= \sum_{i=1}^k \partial \psi_i(x)
	:= \left\{ \sum_{i=1}^k w_i \ \middle|\ w_i \in \partial \psi_i(x), \ i=1,\dots,k \right\}.
	\]
\end{theorem}

\section{Characterization of solution}\label{optimality}

Employing the convex analysis tools introduced in the previous section,  
we now focus on the mathematical characterization of the \emph{generalized $(k,m)$-Heron problem}.    
In this section, we investigate the existence and uniqueness of an optimal solution and derive first-order optimality conditions that describe its structure.  
These results provide the theoretical basis for the algorithmic developments and convergence analysis developed in later sections.

\subsection{Existence of solution}

The generalized $(k,m)$ Heron problem formulated in \eqref{eq:(k,m)HeronProblem} is a convex optimization problem involving the sum of Euclidean norms.  
The objective function is convex and continuous, and each feasible set is assumed to be nonempty, closed, and convex.  
However, convexity alone does not guarantee the existence of a minimizer.  
Existence of a solution requires additional conditions, typically boundedness or coercivity.  

Let $Z = (x_1, \ldots, x_k, y_1, \ldots, y_m) $ denote the vector of all decision variables, and define the feasible region as the Cartesian product
\[
A := S_1 \times S_2 \times \cdots \times S_k \times C_1 \times \cdots \times C_m.
\]
Since each $S_i$ and $C_j$ is nonempty, closed, and convex, the product set $A$ is also nonempty, closed, and convex.  
Thus, problem \eqref{eq:(k,m)HeronProblem} is a constrained convex optimization problem, and any feasible minimizer (if it exists) is globally optimal.

Before establishing the existence of a solution to problem \eqref{eq:(k,m)HeronProblem}, 
we recall a classical principle in convex optimization:  
a continuous function minimized over a closed set attains its minimum whenever the objective is coercive or the feasible region is compact.  

In the present formulation, the feasible region is closed and convex, although the feasible region is not guaranteed compact. However, existence can be ensured under a mild boundedness condition.

\begin{theorem}\label{thm:existence}
	An optimal solution to the generalized $(k,m)$ Heron problem \eqref{eq:(k,m)HeronProblem} exists 
	if at least one of the feasible sets $\{S_1, \ldots, S_k\}$ or the target sets $\{C_1, C_2, \ldots, C_m\}$ is bounded.
\end{theorem}

{\it Proof}
Let 
$ \inf f(Z) = \nu$ 	be the infimum of the objective function over the feasible set $A$.  
Since $f(Z) \geq 0$, $v$ is finite.  

Consider a minimizing sequence 
\[
\{ Z^{(\ell)} \}_{\ell=1}^{\infty} 
= \left\{ \big( x_1^{(\ell)}, \ldots, x_k^{(\ell)}, y_1^{(\ell)}, \ldots, y_m^{(\ell)} \big) \right\},
\]
such that $Z^{(\ell)} \in A$ for all $\ell$ and 
\[
\lim_{\ell \to \infty} f(Z^{(\ell)}) = \nu_{\min}.
\]

Since $f(Z^{(\ell)})$ converges to a finite value, it must be bounded.  
Therefore, there exists some $M > 0$ such that 
\[
f(Z^{(\ell)}) \leq M \qquad \forall~ \text{sufficiently large } \ell.
\]

Assume, without loss of generality, that one of the sets, say $S_p$ for some 
$p \in \{1,\ldots,k\}$, is bounded.  
Because $S_p$ is bounded, there exists a finite constant $B_p$ such that 
\[
\|x_p^{(\ell)}\| \leq B_p \qquad \forall~ x_p^{(\ell)} \in S_p.
\]

From the boundedness of $f(Z^{(\ell)})$, we know that each individual term 
in the sum is bounded: 
\[
\|x_p^{(\ell)} - y_j^{(\ell)}\| < f(Z^{(\ell)}) \leq M, 
\qquad \forall j=1,\ldots,m.
\]

Using the triangle inequality, we have 
\[
\|y_j^{(\ell)}\| 
\leq \|y_j^{(\ell)} - x_p^{(\ell)}\| + \|x_p^{(\ell)}\|.
\]

Since $\|x_p^{(\ell)}\| \leq B_p$, and $\|y_j^{(\ell)} - x_p^{(\ell)}\|$ is bounded, 
it follows that $\|y_j^{(\ell)}\|$ is also bounded for all $j=1,\ldots,m$.
\[
\|y_j^{(\ell)}\| \leq M + B_p
\]

This implies that all sequences $\{y_j^{(\ell)}\}_{\ell=1}^{\infty}$ are bounded.

Now, consider any $x_i^{(\ell)}$ for $i \neq p$.  
For each such $x_i^{(\ell)}$, there is a corresponding term involving $y_j^{(\ell)}$ in the objective function.  

For example, consider $\|x_i^{(\ell)} - y_j^{(\ell)}\|$.  
Since $\|x_i^{(\ell)} - y_j^{(\ell)}\|$ is bounded by $M$ (as part of $f(Z^{(\ell)})$) and 
$\|y_j^{(\ell)}\|$ is now known to be bounded, we can apply the triangle inequality again:  

\[
\|x_i^{(\ell)}\| \leq \|x_i^{(\ell)} - y_j^{(\ell)}\| + \|y_j^{(\ell)}\| 
\leq M + (M + B_p) = 2M + B_p.
\]

This shows that all sequences $\{x_i^{(\ell)}\}_{\ell=1}^{\infty}$ (including $x_p^{(\ell)}$ and those for $i \neq p$) are bounded.  

Consequently, the entire minimizing sequence 
\[
\{Z^{(\ell)}\}_{\ell=1}^{\infty} 
= \big\{(x_1^{(\ell)}, \ldots, x_k^{(\ell)}, y_1^{(\ell)}, \ldots, y_m^{(\ell)}) \big\}_{\ell=1}^{\infty}
\]
is bounded.

By the Bolzano–Weierstrass theorem, a bounded sequence in $A$ must have a subsequence converging to a point in $A$. Let this convergent subsequence be denoted 
$\{Z^{(\ell_r)}\},$ and let its limit be 
\[
Z^* = (x_1^*, \ldots, x_k^*, y_1^*, \ldots, y_m^*).
\]

Since $A$ is closed, $Z^* \in A$. Because the objective function $f(Z)$ is continuous, we also have
\[
f(Z^*) = \lim_{r \to \infty} f(Z^{(\ell_r)}) = \nu_{\min}.
\]

Thus, $Z^*$ is an optimal solution to the problem \eqref{eq:(k,m)HeronProblem}.

Finally, Since the Problem (\ref{eq:(k,m)HeronProblem}) is symmetric with respect to the feasible and target sets.  
If instead of $S_p$ one of the target sets $C_q$ is bounded, the same reasoning applies by interchanging the roles of $S_p$ and $C_q$.  
Hence, an optimal solution to problem~\eqref{eq:(k,m)HeronProblem} exists whenever at least one of the sets 
$S_1, \ldots, S_k, C_1, \ldots, C_m$ is bounded.

\qed

\subsection{Example: Existence of multiple optimal solutions}

The existence theorem established in the previous subsection guarantees that a solution to the generalized $(k,m)$ Heron problem exists under mild boundedness assumptions.  
However, such a solution need not be unique.  
The following symmetric configuration illustrates a case with infinitely many optimal solutions.

\begin{example}
	Consider problem~\eqref{eq:(k,m)HeronProblem} with two feasible sets and two target sets in $\mathbb{R}^2$.  
	The feasible sets are discs centered at $(-6,6)$ and $(6,6)$, each of radius~$2$, the target sets are two smaller discs centered at $(-2,6)$ and $(2,6)$, each of radius~$1$:
	All four sets lie on the horizontal line $y=6$, forming a perfectly symmetric configuration about the $y$-axis.  
	The shortest line segments connecting the boundaries of $S_1$–$C_1$ and $S_2$–$C_2$ coincide with this line.  
	Consequently, $S_1$ and $S_2$ each have unique minimizing points on their rightmost and leftmost boundaries, respectively, while every point along the diameters of $C_1$ and $C_2$ on the line $y=6$ yields the same total distance.  
	Hence, the problem possesses infinitely many optimal solutions.
\end{example}

\begin{figure}[H]
	\centering
	\includegraphics[width=0.8\textwidth]{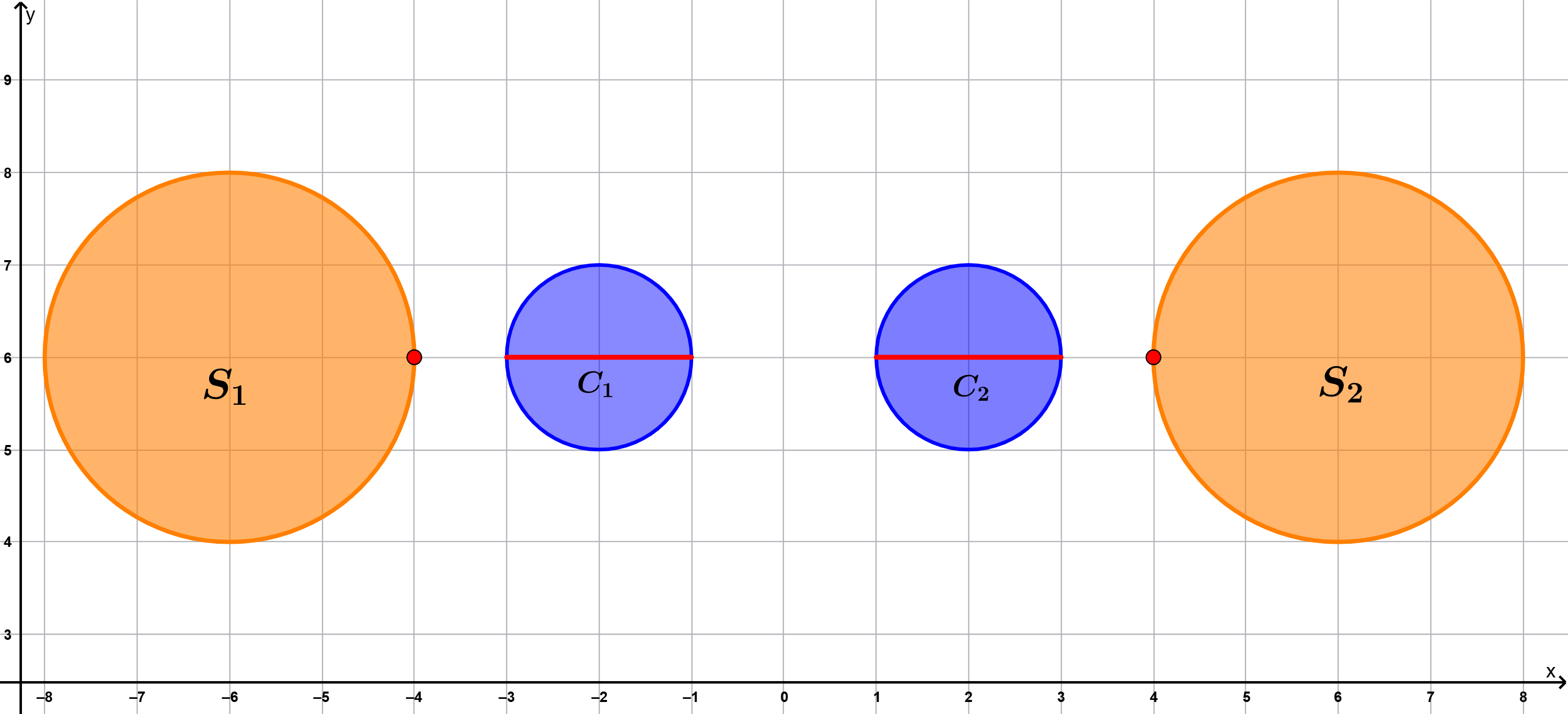}
	\caption{Illustration of multiple optimal solutions to the generalized $(k,m)$-Heron problem.}
	\label{fig:multiple_solution}
\end{figure}

This example shows that, although existence is guaranteed by Theorem~\ref{thm:existence}, uniqueness fails in certain symmetric cases.  
The following subsection establishes additional convexity conditions under which the optimal solution becomes unique.

\subsection{Uniqueness of solution}

While the generalized $(k,m)$-Heron problem may admit multiple solutions in general,  
additional geometric and convexity conditions can ensure uniqueness.  
In this subsection, we develop supporting remarks and lemmas that lead to a sufficient condition for uniqueness of the optimal configuration.

\medskip
Consider the generalized Heron problem  (as discussed in ~\cite{mordukhovich2012convex}):
\[
\min_{x\in S} D(x)
:= \sum_{i=1}^n d(x,C_i)
\quad \text{s.t. } x\in S,
\]
where all $C_i$ and $S$ are nonempty, closed, and convex subsets of $\mathbb{R}^n$, with $C_i\cap S=\emptyset$ for all $i=1,\dots,n$.  
For any $\bar x \in S$, define
\[
a_i(\bar x) := \frac{\bar x - \operatorname{proj}_{C_i}(\bar x)}{d(\bar x,C_i)} \neq 0.
\]
Then, following the result in~\cite[Theorem~3.2, p.~91]{mordukhovich2012convex}, $\bar x$ is an optimal solution if and only if
\begin{equation}\label{Borish-Condition}
	-\sum_{i=1}^n a_i(\bar x) \in N(\bar x;S),
\end{equation}
where $N(\bar x;S)$ denotes the normal cone to $S$ at $\bar x$.

\begin{remark}\label{Boundary_Remarks}
	If each $C_i=\{x_i\}$ is a singleton and all $\{x_i\}_{i=1}^n$ are separated from $S$ by a hyperplane $H$,  
	then all $x_i$ lie on one side of $H$ while $S$ lies on the opposite side.  
	In this case,
	\[
	\sum_{i=1}^n a_i(\bar x) = \sum_{i=1}^n \frac{\bar x -x_i}{\|\bar x-x_i\|} \neq 0.
	\]
	Note that if $\bar x\in\operatorname{int}(S)$, then $N_S(\bar x)=\{0\}$, contradicting the inclusion \eqref{Borish-Condition}.  
	Hence,
	\[
	-\sum_i a_i(\bar x)\in N_S(\bar x)
	\quad \Longleftrightarrow \quad 
	\bar x \in \operatorname{bd}(S).
	\]
	Consequently, if $\bar x\in \operatorname{int}(S)$, then 
	\[
	D(\bar x) > \nu := \min_{x\in S} D(x).
	\]
\end{remark}

\begin{lemma}\label{boundary_lemma}
	Suppose the collections $\{S_i\}_{i=1}^k$ and $\{C_j\}_{j=1}^m$ are nonempty, closed, and convex sets separated by a hyperplane~$H$.  
	If, for a given point 
	\[
	Z = (x_1,\ldots,x_k,y_1,\ldots,y_m),
	\]
	at least one component satisfies $x_i\in \operatorname{int}(S_i)$ (or $y_j\in \operatorname{int}(C_j)$), 
	then there exists another point 
	\[
	\bar Z = (\bar x_1,\ldots,\bar x_k,\bar y_1,\ldots,\bar y_m)
	\quad \text{such that} \quad 
	F(\bar Z) < F(Z).
	\]
\end{lemma}

{\it Proof}
By the preceding remark~\ref{Boundary_Remarks}, each minimizer must lie on the boundary of its feasible set.  
If some $x_i\in\operatorname{int}(S_i)$, one can select $\bar x_i\in\operatorname{bd}(S_i)$ such that the total distance strictly decreases, i.e., $F(\bar Z) < F(Z)$.  
\qed

\begin{theorem}[Uniqueness of Solution]
	\label{thm:uniqueness}
	Suppose that the collections $\{S_i\}_{i=1}^k$ and $\{C_j\}_{j=1}^m$ are strictly convex and separated by a hyperplane~$H$.  
	Then the generalized $(k,m)$-Heron problem~\eqref{eq:(k,m)HeronProblem} admits a unique optimal solution.
\end{theorem}

{\it Proof}
Let 
$Z = (x_1,\ldots,x_k,y_1,\ldots,y_m)$,   $W = (p_1,\ldots,p_k,q_1,\ldots,q_m)$ be two distinct optimal solutions such that $F(Z)=F(W)=\nu$.  
Then $x_i\neq p_i$ for some~$i$.  
Since each $S_i$ is strictly convex, 
\[
\frac{x_i+p_i}{2} \in \operatorname{int}(S_i),
\]
and by the previous Lemma \ref{boundary_lemma}, 
\[
F\!\left(\frac{Z+W}{2}\right) > F(Z)=\nu.
\]
However, since $F$ is convex,
\[
F\!\left(\frac{Z+W}{2}\right) \le \frac{F(Z)+F(W)}{2} = \nu,
\]
which is a contradiction.  
Therefore, the generalized $(k,m)$-Heron problem admits a unique optimal solution.
\qed

\subsection{Optimality conditions}

The preceding examples establish that while an optimal solution to the generalized $(k,m)$ Heron problem always exists under mild assumptions, it may not be unique.  
To further understand the nature and structure of optimal solutions, we now derive the necessary and sufficient conditions that characterize them.  
These conditions describe the geometric balance between the feasible and target sets at optimality and serve as the foundation for the algorithmic methods developed later in this paper.

We begin by reformulating the constrained optimization problem into an equivalent unconstrained formulation and then applying subgradient calculus to obtain the first-order optimality conditions.

We can rewrite problem~\eqref{eq:(k,m)HeronProblem} equivalently as
\begin{equation}\label{eq:UnconstrainedObjective}
	\min_{Z \in \mathbb{R}^{n(k+m)}} \Big\{ f(Z) + \delta_A(Z) \Big\},	
\end{equation}
where $\delta_A(Z)$ denotes the indicator function of the feasible set $A$.  
This reformulation allows us to express the optimality conditions compactly in terms of subgradients and normal cones.

\begin{theorem}[First-Order Optimality Conditions]
	\label{thm:optimality}
	Let $S_1, \ldots, S_k$ and $C_1, \ldots, C_m$ be nonempty, closed, and convex subsets of $\mathbb{R}^n$ and $S_i \cap C_j = \varnothing$ for all $i, j.$  
	Let $Z^* = (x_1^*, \ldots, x_k^*,\, y_1^*, \ldots, y_m^*)$ be a feasible point, where $x_i^* \in S_i,$ $i = 1, \ldots, k$ and $y_j^* \in C_j,$ $j = 1, \ldots, m.$  
	Then $Z^*$ is an optimal solution of problem~\eqref{eq:(k,m)HeronProblem} if and only if there exist normal vectors $n_{S_i} \in N_{S_i}(x_i^*)$ and $n_{C_j} \in N_{C_j}(y_j^*)$ satisfying the following conditions.

	\medskip
	\noindent\textbf{(O1)} For each $i=1,\ldots,k$,
	\[
	\sum_{j=1}^{m} \frac{x_i^* - y_j^*}{\Vert x_i^* - y_j^* \Vert} + n_{S_i} = 0,
	\]
	
	\medskip
	\noindent\textbf{(O2)} For each $j=1,\ldots,m$,
	\[
	\sum_{i=1}^{k} \dfrac{y_j^*-x_i^*}{\|y_j^*-x_i^*\|} + n_{C_j} = 0,
	\]
	
	\medskip
	\noindent\textbf{(O3)} 	Summing (O1) over $i=1,\dots,k$ and (O2) over $j=1,\dots,m$ yields the global equilibrium condition
	\[
	\sum_{i=1}^k n_{S_i} + \sum_{j=1}^m n_{C_j} = 0,
	\]
	which expresses the balance of normal vectors across all feasible and target sets.
	
\end{theorem}
{\it Proof} 
We reformulate this generalized $(k,m)$ Heron problem~\eqref{eq:(k,m)HeronProblem} into an equivalent unconstrained problem by using the indicator function. Let $Z = (x_1, x_2, \ldots, x_k, y_1, y_2, \ldots, y_m)$ $\in \mathbb{R}^{n(k+m)}$ be the decision vector of all variables. The feasible region is $A := S_1 \times S_2 \times \cdots \times S_k \times C_1 \times \cdots \times C_m,$ which is nonempty, closed, and convex. Introducing the indicator function $\delta_A$, we can represent \eqref{eq:(k,m)HeronProblem} equivalently as

\[
\min_{Z \in \mathbb{R}^{n(k+m)}} \; \phi(Z) := F(Z) + \delta_A(Z).
\]

For a convex function $\phi(Z)$, a point $Z^*$ is a global minimizer if and only if 
the zero vector belongs to its subdifferential at $Z^*$:
\begin{equation}\label{zero_belongs_Sub}
	0 \in \partial \phi(Z^*) 
	= \partial F(Z^*) + \partial \delta_A(Z^*)
\end{equation}

Recall that the subdifferential of the indicator function $\delta_A(Z)$ is the 
normal cone to the feasible set $A$. Moreover, by the product structure of $A$, 
we can apply the Proposition \ref{Normalcone_Product} to write
\[
\partial \delta_A(Z^*) = N_A(Z^*) 
= N_{S_1}(x_1^*) \times \cdots \times N_{S_k}(x_k^*) \times N_{C_1}(y_1^*) \times \cdots \times N_{C_m}(y_m^*).
\]
Thus, any $n \in N_A(Z^*)$ has the vector form
$n = (n_{S_1}^*, \ldots, n_{S_k}^*, \, n_{C_1}^*, \ldots, n_{C_m}^*),$ where $n_{S_i}^* \in N_{S_i}(x_i^*)$ and $n_{C_j}^* \in N_{C_j}(y_j^*)$.

Since
\[
F(Z) = \sum_{i=1}^k \sum_{j=1}^m \|x_i - y_j\|.
\]

For each $i = 1, \ldots, k$, \\
From the Affine rule Theorem \ref{affinetheorem}, we have
\[
\partial_{x_i} F(Z^*) = \sum_{j=1}^m \frac{x_i^* - y_j^*}{\|x_i^* - y_j^*\|}.
\]

Similarly, for each $j = 1, \ldots, m$,
\[
\partial_{y_j} F(Z^*) = \sum_{i=1}^k \frac{y_j^* - x_i^*}{\|y_j^* - x_i^*\|}.
\]

Thus,
\[
\partial F(Z^*) = \big( \partial_{x_1} F(Z^*), \ldots, \partial_{x_k} F(Z^*), 
\partial_{y_1} F(Z^*), \ldots, \partial_{y_m} F(Z^*) \big).
\]

\[
\partial F(Z^*) =
\Bigg(
\sum_{j=1}^m \frac{x_1^* - y_j^*}{\|x_1^* - y_j^*\|}, \;
\ldots, \;
\sum_{j=1}^m \frac{x_k^* - y_j^*}{\|x_k^* - y_j^*\|}, \;
\sum_{i=1}^k \frac{y_1^* - x_i^*}{\|x_i^* - y_1^*\|}, \;
\ldots, \;
\sum_{i=1}^k \frac{y_m^* - x_i^*}{\|x_i^* - y_m^*\|}
\Bigg).
\]

Let $\mathbf{0}$ be represented as a tuple of zero vectors of dimension $n$:  
\[
\mathbf{0} = (\mathbf{0}_n, \ldots, \mathbf{0}_n, \mathbf{0}_n, \ldots, \mathbf{0}_n).
\]

Thus the condition \eqref{zero_belongs_Sub} becomes,

\begin{equation*}
	\begin{aligned}
		(\mathbf{0}_n, \ldots, \mathbf{0}_n) \in &
		\Big(
		\sum_{j=1}^m \tfrac{x_1^* - y_j^*}{\|x_1^* - y_j^*\|}, \ldots,
		\sum_{j=1}^m \tfrac{x_k^* - y_j^*}{\|x_k^* - y_j^*\|}, \\[4pt]
		&\quad
		\sum_{i=1}^k \tfrac{y_1^* - x_i^*}{\|x_i^* - y_1^*\|}, \ldots,
		\sum_{i=1}^k \tfrac{y_m^* - x_i^*}{\|x_i^* - y_m^*\|}
		\Big)
		+ (n_{S_1}, \ldots, n_{S_k}, n_{C_1}, \ldots, n_{C_m}).
	\end{aligned}
\end{equation*}

For this inclusion to hold, each component must sum to the zero vector, which gives:

For $i = 1, \ldots, k$
\[
\mathbf{0} = \sum_{j=1}^m \frac{x_i^* - y_j^*}{\|x_i^* - y_j^*\|} + n_{S_i},
\tag{O1}
\]

For $j = 1, \ldots, m$
\[
\mathbf{0} = \sum_{i=1}^k \frac{y_j^* - x_i^*}{\|x_i^* - y_j^*\|} + n_{C_j}.
\tag{O2}
\]

This can be rearranged as:
\[
n_{S_i} = - \sum_{j=1}^m \frac{x_i^* - y_j^*}{\|x_i^* - y_j^*\|}
= \sum_{j=1}^m \frac{y_j^* - x_i^*}{\|x_i^* - y_j^*\|},
\]
\[
n_{C_j} = - \sum_{i=1}^k \frac{y_j^* - x_i^*}{\|x_i^* - y_j^*\|}
= \sum_{i=1}^k \frac{x_i^* - y_j^*}{\|x_i^* - y_j^*\|}.
\]

Summing \((O1)\) over $k$ and \((O2)\) over $m$ gives
\[
\sum_{i=1}^k n_{S_i}^* + \sum_{j=1}^m n_{C_j}^* = 0.
\tag{O3}
\]
which completes the proof.
\qed

\section{Numerical Algorithm}\label{numerical-algorithm}

The optimality conditions derived in the previous section provide a theoretical characterization of the solution structure for the generalized $(k,m)$ Heron problem.  
To approximate such optimal points, we now design an iterative procedure based on subgradient information and projection operators.  
In the following, we present the \emph{Projected Subgradient Algorithm (PSA)}, for our proposed problem \eqref{eq:(k,m)HeronProblem}, and establish its convergence and  subsequently demonstrate its performance through illustrative numerical examples.

Before presenting the proof of convergence, recall that in the analysis of subgradient-based methods it is standard to assume that the norm of the subgradients of the objective function is uniformly bounded (see, for example,~\cite[Assumption 3.2.1, p. 153]{bertsekas2015convex}).  
Formally, there exists a constant $G>0$ such that 
$\|g_t\| \le G$ for all iterations $t$.  
This assumption guarantees the stability of the projected subgradient updates and 
is automatically satisfied when the objective function is Lipschitz continuous with Lipschitz constant~$G$.

\begin{remark}\label{rem:subgrad-bound}
	For our generalized $(k,m)$-Heron problem \eqref{eq:(k,m)HeronProblem}, this boundedness can be verified explicitly. Consider our objective function \eqref{eq:(k,m)HeronProblem}, 
	$F(Z)=\sum_{i=1}^{k}\sum_{j=1}^{m}\|x_i-y_j\|,\;
	Z=(x_1,\ldots,x_k,\allowbreak y_1,\ldots,y_m)$.
	The corresponding subgradients are at each iteration $t$, the vector $g^{(t)}$ is a subgradient of the convex function $F$ at $Z^{(t)}$, given by
	$g^{(t)} = (g_{x_1}^{(t)}, \ldots, g_{x_k}^{(t)}, g_{y_1}^{(t)}, \ldots, g_{y_m}^{(t)}),$
	where each component is computed as
	$g_{x_i}^{(t)} = \sum_{j=1}^m u_{ij}^{(t)},$  $g_{y_j}^{(t)} = \sum_{i=1}^k v_{ji}^{(t)},$  with
	\[
	u_{ij}^{(t)} =
	\begin{cases}
		\dfrac{x_i^{(t)} - y_j^{(t)}}{\|x_i^{(t)} - y_j^{(t)}\|}, & \text{if } x_i^{(t)} \ne y_j^{(t)}, \\[0.8em]
		0, & \text{if } x_i^{(t)} = y_j^{(t)},
	\end{cases}
	\qquad
	v_{ji}^{(t)} =
	\begin{cases}
		\dfrac{y_j^{(t)} - x_i^{(t)}}{\|y_j^{(t)} - x_i^{(t)}\|}, & \text{if } y_j^{(t)} \ne x_i^{(t)}, \\[0.8em]
		0, & \text{if } y_j^{(t)} = x_i^{(t)}.
	\end{cases}
	\]
	
	Consequently, the subgradient norm is bounded as
	
	\[
	\|g\|^2
	=\sum_{i=1}^{k}\|g^{(x_i)}\|^2+\sum_{j=1}^{m}\|g^{(y_j)}\|^2
	\le k\,m^2+m\,k^2
	=km(m+k),
	\]
	and thus every subgradient of $F$ obeys
	\[
	\|g\|\le G:=\sqrt{km(m+k)}.
	\]
	This explicit constant $G$ plays the same role as the bound $c$ in the classical
	subgradient analysis~\cite[Proposition 3.2.1, p.~153 ]{bertsekas2015convex}, ensuring the
	Projected Subgradient Algorithm convergence properties established below.
	
\end{remark}

\begin{theorem}[Convergence of the PSA]\label{thm:PSA_convergence}
	Let $S_1,\ldots,S_k$ and $C_1,\ldots,C_m$ be nonempty, closed, and convex subsets of $\mathbb{R}^n$, and consider the generalized $(k,m)$-Heron problem
	\[
	\min_{Z\in A} \; F(Z)=\sum_{i=1}^{k}\sum_{j=1}^{m}\|x_i-y_j\|,
	\qquad 
	A = S_1\times\cdots\times S_k\times C_1\times\cdots\times C_m.
	\]
	Assume that at least one of the sets among $\{S_i\}$ and $\{C_j\}$ is bounded, guaranteeing the existence of an optimal solution $Z^*\in A$ with optimal value $F^*:=F(Z^*)$.  
	
	Let $\{Z_t\}$ be the sequence generated by the \emph{Projected Subgradient Algorithm}:
	\begin{equation}\label{eq:PSA_update}
		Z_{t+1} = \proj_A(Z_t - \alpha_t g_t),
	\end{equation}
	where $g_t \in \partial F(Z_t)$, $\proj_A$ denotes the Euclidean projection onto $A$, and the step-size sequence $\{\alpha_t\}$ satisfies
	\[
	\alpha_t > 0, \qquad 
	\sum_{t=0}^\infty \alpha_t = \infty, \qquad 
	\sum_{t=0}^\infty \alpha_t^2 < \infty.
	\]
	Then the best-iterate sequence of objective values
	$F_{\mathrm{best}}^{(N)} := \min_{0\le t<N} F(Z_t)$
	converges to the optimal value, i.e.,$\lim_{N\to\infty} F_{\mathrm{best}}^{(N)} = F^*$, and $Z_t$  converges to an optimal solution.

\end{theorem}
{\it Proof} 
The existence of an optimal solution $Z^*\in A$ follows from Theorem~\ref{thm:existence}.  
By the nonexpansiveness of Euclidean projections onto closed convex sets \cite[Proposition 3.2.1, p.~149]{bertsekas2015convex},
\begin{equation}\label{convergence-equation}
	\|Z_{t+1}-Z^*\|^2
	=\|\proj_A(Z_t-\alpha_t g_t)-\proj_A(Z^*)\|^2
	\le \|Z_t-\alpha_t g_t - Z^*\|^2
\end{equation}
Expanding and applying the subgradient inequality $F(Z_t)-F(Z^*) \le g_t^\top(Z_t-Z^*)$ gives
\begin{equation}\label{eq:basic-rec}
	\|Z_{t+1}-Z^*\|^2
	\le \|Z_t-Z^*\|^2 - 2\alpha_t\bigl[F(Z_t)-F^*\bigr] + \alpha_t^2\|g_t\|^2.
\end{equation}
Summing~\eqref{eq:basic-rec} from $t=0$ to $N-1$ and using $\|g_t\|\le G$ yields
\begin{equation}\label{eq:sum-bound}
	\sum_{t=0}^{N-1}\alpha_t\bigl[F(Z_t)-F^*\bigr]
	\le \tfrac{1}{2}\|Z_0-Z^*\|^2 + \tfrac{G^2}{2}\sum_{t=0}^{N-1}\alpha_t^2.
\end{equation}
Let $F_{\mathrm{best}}^{(N)} = \min_{0\le t<N} F(Z_t)$.  
Since $F(Z_t)-F^*\ge F_{\mathrm{best}}^{(N)}-F^*$ for all $t \leq N$, \eqref{eq:sum-bound} implies
\[
0 \leq 	F_{\mathrm{best}}^{(N)}-F^*
\le \frac{\|Z_0-Z^*\|^2 + G^2\sum_{t=0}^{N-1}\alpha_t^2}
{2\sum_{t=0}^{N-1}\alpha_t}.
\]
Under the step-size assumptions \cite[Proposition 3.2.6, p.~157]{bertsekas2015convex},  $\sum_t \alpha_t = \infty$ and $\sum_t \alpha_t^2 < \infty$, the numerator remains bounded while the denominator diverges, therefore
$\lim_{N\to\infty} \bigl(F_{\mathrm{best}}^{(N)} - F^*\bigr) = 0,$
that is $F_{\mathrm{best}}^{(N)}\to F^*$.

Using Proposition A.2.4 \cite[See p.~462]{bertsekas2015convex} or Proposition A.2.6 \cite[See p.~465]{bertsekas2015convex}, for our settings, in particular identifying $\beta_k ,\gamma_k$  and  $\delta_k$ in \eqref{convergence-equation}, and set $\phi(Z_t, Z^*)=F(Z_t)-F(Z^*)$, we see that $Z_t$ converges to some optimal point.



\qed
The convergence of the Projected Subgradient Algorithm (PSA) was established in Theorem~\ref{thm:PSA_convergence}. 
We now outline its practical implementation for solving the generalized $(k,m)$-Heron problem. 
In all computations, the diminishing step-size rule 
\begin{equation}
	\alpha_t \downarrow 0,\; \sum_t \alpha_t = \infty,\; \sum_t \alpha_t^2 < \infty
\end{equation}
is adopted, as it satisfies the standard convergence conditions for subgradient methods. 

The complete iterative procedure of the PSA is summarized in Algorithm~\ref{alg:PSA_km}.
\begin{algorithm}[htbp]
	\caption{Projected Subgradient Algorithm (PSA) for the Generalized $(k,m)$-Heron Problem}
	\label{alg:PSA_km}
	\KwIn{
		Closed and convex sets $S_i \subset \mathbb{R}^n$ for $i=1,\dots,k$ and $C_j \subset \mathbb{R}^n$ for $j=1,\dots,m$;\\
		Initial feasible point $Z^{(0)} = (x_1^{(0)},\dots,x_k^{(0)},y_1^{(0)},\dots,y_m^{(0)})$ with $x_i^{(0)}\in S_i$ and $y_j^{(0)}\in C_j$;\\
		Step-size sequence $\{\alpha_t\}_{t\ge0}$ satisfying the diminishing rule 
		$\alpha_t \downarrow 0$, 
		$\sum_t \alpha_t = \infty$, 
		$\sum_t \alpha_t^2 < \infty$;\\
		Tolerance $\varepsilon > 0$.
	}
	
	\KwOut{Approximate optimal solution $Z^*=(x_1^*,\dots,x_k^*,y_1^*,\dots,y_m^*)$ and objective value $F^*$.}
	
	\BlankLine
	\For{$t = 0,1,2,\dots$}{
		\tcp{Compute subgradients with respect to $x_i$}
		\For{$i = 1$ \KwTo $k$}{
			\For{$j = 1$ \KwTo $m$}{
				$u_{ij}^{(t)} \leftarrow
				\begin{cases}
					\dfrac{x_i^{(t)} - y_j^{(t)}}{\|x_i^{(t)} - y_j^{(t)}\|}, & \text{if } x_i^{(t)} \ne y_j^{(t)},\\[0.4em]
					0, & \text{otherwise.}
				\end{cases}$\;
			}
			$g_{x_i}^{(t)} \leftarrow \displaystyle\sum_{j=1}^{m} u_{ij}^{(t)}$\;
		}
		
		\tcp{Compute subgradients with respect to $y_j$}
		\For{$j = 1$ \KwTo $m$}{
			\For{$i = 1$ \KwTo $k$}{
				$v_{ji}^{(t)} \leftarrow
				\begin{cases}
					\dfrac{y_j^{(t)} - x_i^{(t)}}{\|y_j^{(t)} - x_i^{(t)}\|}, & \text{if } y_j^{(t)} \ne x_i^{(t)},\\[0.4em]
					0, & \text{otherwise.}
				\end{cases}$\;
			}
			$g_{y_j}^{(t)} \leftarrow \displaystyle\sum_{i=1}^{k} v_{ji}^{(t)}$\;
		}

		$g^{(t)} \leftarrow (g_{x_1}^{(t)},\dots,g_{x_k}^{(t)},g_{y_1}^{(t)},\dots,g_{y_m}^{(t)})$\;
		$\tilde{z}^{(t+1)} \leftarrow z^{(t)} - \alpha_t g^{(t)}$\;
		
		\tcp{Project onto the feasible set $A$}  
		$Z^{(t+1)} \leftarrow \proj_A(\tilde{z}^{(t+1)})$\;
		
		\tcp{Evaluate objective and check stopping criterion}
		$F^{(t+1)} \leftarrow \displaystyle\sum_{i=1}^{k}\sum_{j=1}^{m} \|x_i^{(t+1)} - y_j^{(t+1)}\|$\;
		
		\If{$|F^{(t+1)} - F^{(t)}| < \varepsilon$}{
			\textbf{stop and return} $Z^{(t+1)}, F^{(t+1)}$\;
		}
	}
\end{algorithm}
\clearpage

\section{Numerical Illustration}\label{sec:numerical-example}

To validate and illustrate the theoretical findings developed in the previous section, 
we now turn to numerical experiments based on the proposed Projected Subgradient Method (PSM).  
These examples not only confirm the convergence guarantees established earlier 
but also reveal the algorithm’s practical efficiency and geometric insight. We present two illustrative examples of the generalized $(k,m)$-Heron problem.
The first example deals with a two-dimensional configuration in $\mathbb{R}^2$,
while the second extends the study to a three-dimensional setting in $\mathbb{R}^3$,
demonstrating the method’s scalability and its ability to handle more complex geometric structures.
\paragraph{Experimental Setup.}
All experiments were implemented in Python~3.11 on Google~Colab using an Intel\textsuperscript{\textregistered} Xeon\textsuperscript{\textregistered} CPU @~2.20\,GHz with 12.7\,GB of RAM. The \emph{Projected Subgradient Algorithm (PSA)} employed the step-size rule $\alpha_t = 1/t$, which satisfies the diminishing conditions $\alpha_t \downarrow 0$, $\sum_t \alpha_t = \infty$, and $\sum_t \alpha_t^2 < \infty$. The stopping criterion was set as $|F(x^{(k+1)}) - F(x^{(k)})| < 10^{-15}$, where $F(x)$ denotes the objective function in~\eqref{eq:(k,m)HeronProblem}.
\begin{example}[Two-Dimensional $(4,3)$-Heron Problem]
	\label{ex:2D_Heron}
	We consider the generalized $(k,m)$-Heron problem in $\mathbb{R}^2$ with $k=4$ and $m=3$. 
	The feasible sets $\{S_i\}_{i=1}^4$ are unit discs centered at $(8,9)$, $(2,9)$, $(-2,12)$, and $(-7,8)$, while the target sets $\{C_j\}_{j=1}^3$ are squares of half-side length~$1$ centered at $(4,2)$, $(6,12)$, and $(-3,6)$. 
	The initial points, chosen on the respective set boundaries (not at the centers), are 
	$x_1^{(0)}=(9,5)$, $x_2^{(0)}=(2,10)$, $x_3^{(0)}=(-2,13)$, $x_4^{(0)}=(-8,8)$, 
	and $y_1^{(0)}=(5,1)$, $y_2^{(0)}=(7,13)$, $y_3^{(0)}=(-4,5)$. 
	The \emph{Projected Subgradient Algorithm (PSA)} in Algorithm~\ref{alg:PSA_km} is applied with a diminishing step-size rule $\alpha_t = 1/t$ and tolerance $\varepsilon = 10^{-15}$.
	
	\smallskip
	The algorithm converged after approximately $1.69\times10^5$ iterations, yielding the optimal objective value $F^* = 79.113613$. 
	The optimal points are 
	$x_1^*=(7.0399,5.2796)$, $x_2^*=(1.9216,8.0031)$, $x_3^*=(-1.4238,11.1827)$, $x_4^*=(-6.0103,7.8565)$, 
	and $y_1^*=(3.0000,3.0000)$, $y_2^*=(5.0000,11.0000)$, $y_3^*=(-2.0000,7.0000)$. 
	All optimal points lie on their respective set boundaries, satisfying the theoretical optimality conditions. 
	The convergence history (Table~\ref{tab:conv_summary_ex1}) indicates rapid early reduction in $F^{(t)}$ and stabilization near $F^*$, confirming the numerical accuracy and theoretical convergence of the PSA. 
	The pairwise Euclidean distances between $(x_i^*,y_j^*)$ are reported in Table~\ref{tab:distance_matrix_example1}, where the smallest values correspond to $(x_2^*,y_3^*)$ and $(x_4^*,y_3^*)$, reflecting their geometric proximity in the optimal configuration shown in Figure~\ref{fig:example1}.
	
	\begin{figure}[H]
		\centering
		\includegraphics[width=0.9\textwidth]{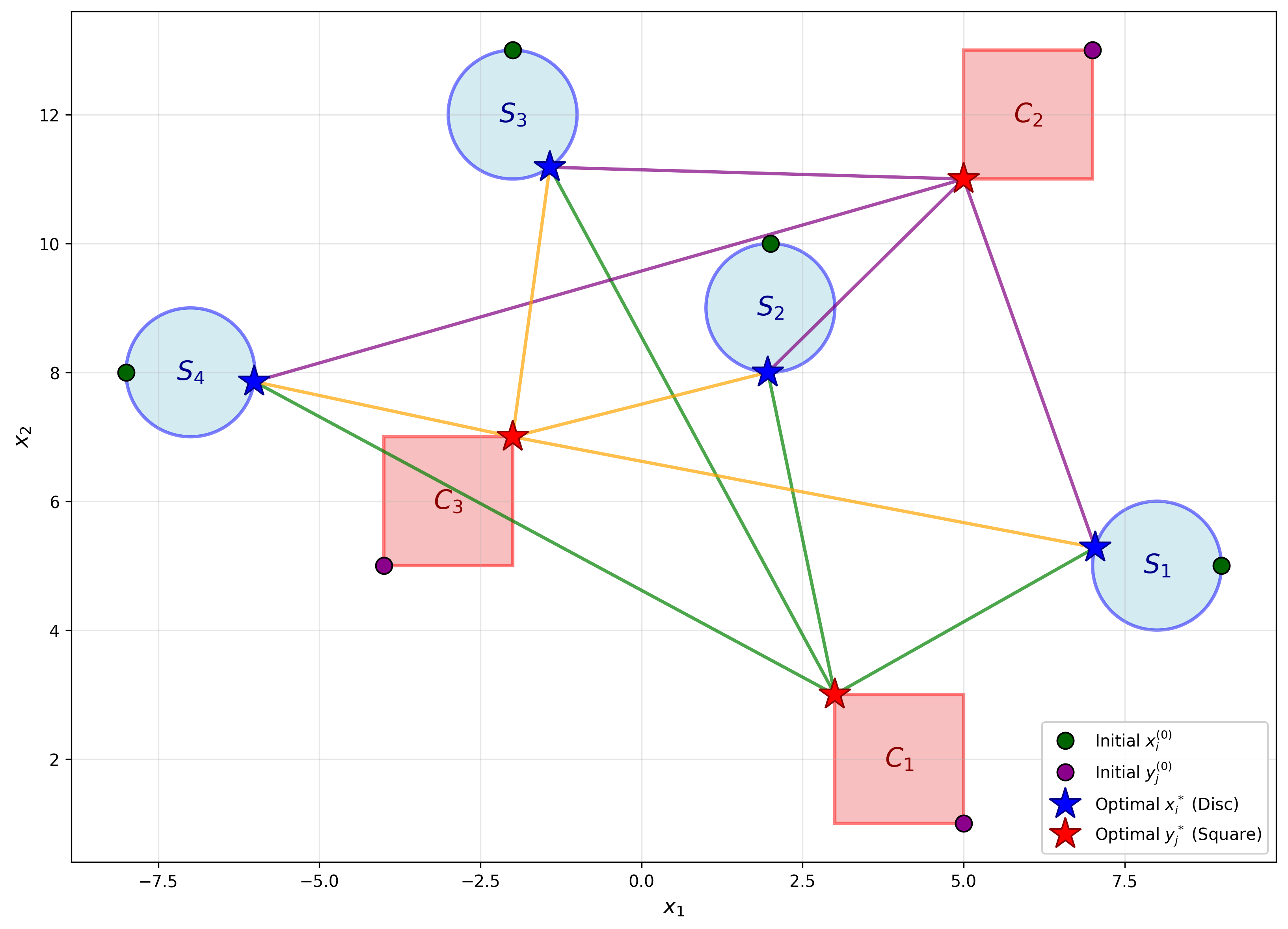}
		\caption{Illustration of Example \ref{ex:2D_Heron}}
		\label{fig:example1}
	\end{figure}
\end{example}

\begin{table}[H]
	\centering
	\caption{Convergence behavior of the PSA for the $(4,3)$-Heron problem in $\mathbb{R}^2$.}
	\label{tab:conv_summary_ex1}
	\begin{tabular}{rcc}
		\toprule
		\textbf{Iteration} & \textbf{Objective Value} $F^{(t)}$ & \textbf{Change} $|\Delta F|$ \\
		\midrule
		1       & 82.117842 & $3.4628\times10^{+1}$ \\
		1{,}000 & 79.113630 & $2.9860\times10^{-8}$ \\
		10{,}000 & 79.113613 & $5.5607\times10^{-11}$ \\
		50{,}000 & 79.113613 & $6.8212\times10^{-13}$ \\
		100{,}000 & 79.113613 & $1.1369\times10^{-13}$ \\
		169{,}449 & 79.113613 & $0.0000$ \\
		\bottomrule
	\end{tabular}
\end{table}

\begin{table}[H]
	\centering
	\caption{Pairwise Euclidean distances $\|x_i^* - y_j^*\|$ between optimal feasible points $x_i^*$ and target points $y_j^*$ for the $(4,3)$-Heron problem in $\mathbb{R}^2$.}
	\label{tab:distance_matrix_example1}
	\begin{tabular}{c|ccc}
		\toprule
		\textbf{Feasible Points} & \textbf{$y_1^*$} & \textbf{$y_2^*$} & \textbf{$y_3^*$} \\
		\midrule
		$x_1^*$ & 4.6386 & 6.0733 & 9.2021 \\
		$x_2^*$ & 5.1180 & 4.2963 & 4.0478 \\
		$x_3^*$ & 9.3020 & 6.4264 & 4.2222 \\
		$x_4^*$ & 10.2358 & 11.4503 & 4.1008 \\
		\bottomrule
	\end{tabular}
\end{table}

\begin{example}[Three-Dimensional $(3,2)$-Heron Problem]
	\label{ex:psa_3d}
	We now consider a three-dimensional instance of the generalized $(k,m)$-Heron problem in~$\mathbb{R}^3$ with $k=3$ and $m=2$.  
	The feasible sets $\{S_i\}_{i=1}^3$ are spheres of radius~$1$ centered at $(-3,1,2)$, $(1,4,4)$, and $(4,1,2)$, 
	while the target sets $\{C_j\}_{j=1}^2$ are cubes with half-side length~$1$ centered at $(-3,-1,-2)$ and $(3,-3,-2)$.  
	The initial points are chosen within the corresponding sets (but not at their centers) as 
	$x_1^{(0)}=(-4,1,2)$, $x_2^{(0)}=(-1,4,5)$, $x_3^{(0)}=(5,1,2)$, 
	and $y_1^{(0)}=(-4,0,-3)$, $y_2^{(0)}=(4,-2,-3)$.  
	The \emph{Projected Subgradient Algorithm (PSA)} from Algorithm~\ref{alg:PSA_km} 
	is applied with a diminishing step-size rule $\alpha_t = 1/t$ and tolerance $\varepsilon = 10^{-15}$.
	
	\smallskip
	The algorithm converged after $1.29\times10^3$ iterations, attaining the optimal objective value $F^* = 30.691348$.  
	The optimal feasible and target points are 
	$x_1^*=(-2.4585,0.6055,1.2576)$, 
	$x_2^*=(0.8422,3.3061,3.2974)$, 
	$x_3^*=(3.3092,0.5701,1.4186)$, 
	and 
	$y_1^*=(-2.0000,0.0000,-1.0000)$, 
	$y_2^*=(2.0000,-2.0000,-1.0000)$.

	All optimal points lie on their respective set boundaries and satisfy the first-order optimality conditions established earlier.  
	The convergence pattern, summarized in Table~\ref{tab:conv_summary_ex2}, shows rapid decay in the objective within the first few iterations, confirming the numerical stability and efficiency of the PSA.  
	The pairwise distances between the optimal points are reported in Table~\ref{tab:distance_matrix_example2}, 
	illustrating that $(x_1^*,y_1^*)$ and $(x_3^*,y_2^*)$ form the shortest connecting pairs, consistent with the geometric configuration shown in Figure~\ref{fig:example2}.
	
	\begin{figure}[H]
		\centering
		\includegraphics[width=1.0\textwidth]{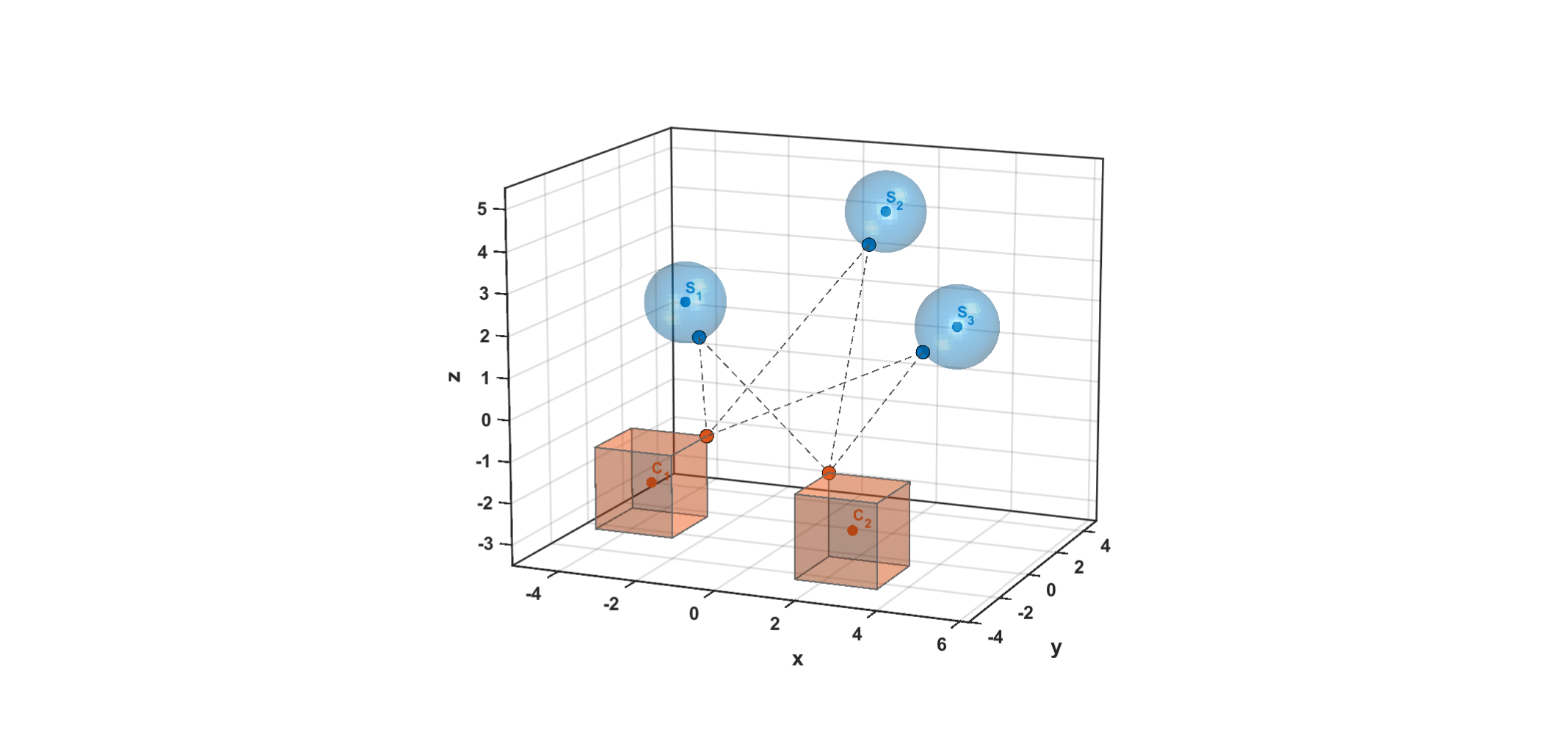}
		\caption{Illustration of Example \ref{ex:psa_3d}}
		\label{fig:example2}
	\end{figure}
	
\end{example}

\begin{table}[H]
	\centering
	\caption{Convergence behavior of the PSA for the $(3,2)$-Heron problem in $\mathbb{R}^3$.}
	\label{tab:conv_summary_ex2}
	\begin{tabular}{rcc}
		\toprule
		\textbf{Iteration} & \textbf{Objective Value} $F^{(t)}$ & \textbf{Change} $|\Delta F|$ \\
		\midrule
		1       & 34.161552 & $1.7711\times10^{+1}$ \\
		10      & 30.693386 & $1.0413\times10^{-3}$ \\
		100     & 30.691348 & $5.7477\times10^{-9}$ \\
		500     & 30.691348 & $1.0303\times10^{-12}$ \\
		1{,}000 & 30.691348 & $2.1316\times10^{-14}$ \\
		1{,}289 & 30.691348 & $0.0000$ \\
		\bottomrule
	\end{tabular}
\end{table}

\begin{table}[H]
	\centering
	\caption{Pairwise Euclidean distances $\|x_i^* - y_j^*\|$ between optimal feasible points $x_i^*$ and target points $y_j^*$ for the $(3,2)$-Heron problem in $\mathbb{R}^3$.}
	\label{tab:distance_matrix_example2}
	\begin{tabular}{c|cc}
		\toprule
		\textbf{Feasible Points} & \textbf{$y_1^*$} & \textbf{$y_2^*$} \\
		\midrule
		$x_1^*$ & 2.3819 & 5.6359 \\
		$x_2^*$ & 6.1218 & 6.9255 \\
		$x_3^*$ & 5.8620 & 3.7642 \\
		\bottomrule
	\end{tabular}
\end{table}

\section{Conclusion}
This study introduced and thoroughly investigated the \emph{generalized $(k,m)$-Heron problem}, extending the classical and generalized Heron formulations to configurations involving multiple feasible and multiple target convex sets. using modern convex analysis tools, we established fundamental theoretical results, including the existence and uniqueness of optimal solutions and precise first-order optimality conditions derived from subdifferential and normal cone analysis. Based on these theoretical foundations, we proposed a \emph{Projected Subgradient Algorithm (PSA)} to numerically solve the problem and rigorously proved its convergence under a standard diminishing step-size rule. Numerical experiments in both $\mathbb{R}^2$ and $\mathbb{R}^3$ verified the analytical findings and highlighted the algorithm’s robustness, efficiency, and geometric accuracy. The presented framework not only unifies several distance-based geometric optimization models but also opens new avenues for applications in areas such as multi-agent localization, geometric design, and network optimization.   

\section*{Author Contributions}
All authors contributed equally to this work.

\section*{Conflict of Interest}
The authors declare that there is no conflict of interest.

\section*{Funding}
This research received no external funding.

\bibliographystyle{plain}
\bibliography{km-Heron-Ref}

\end{document}